\documentclass[aop,preprint]{imsart}

\RequirePackage[OT1]{fontenc}
\RequirePackage{amsthm,amsmath}
\usepackage{amsfonts}
\usepackage{amsmath}
\RequirePackage[numbers]{natbib}
\RequirePackage[colorlinks,citecolor=blue,urlcolor=blue]{hyperref}

\arxiv{arXiv:0000.0000}

\startlocaldefs

\newtheorem {Proposition}{Proposition}[section]
\newtheorem {Lemma}[Proposition] {Lemma}
\newtheorem {Theorem}[Proposition]{Theorem}
\newtheorem {Corollary}[Proposition]{Corollary}
\newtheorem {Remark}[Proposition]{Remark}

\endlocaldefs

\begin{document}

\begin{frontmatter}
\title{Central Limit Theorems for empirical transportation cost in general dimension}

\runtitle{CLTS for empirical transportation costs}
%\thankstext{T1}{Footnote to the title with the ``thankstext'' command.}

\begin{aug}
\author{\fnms{Eustasio} \snm{ del Barrio}\ead[label=e1]{tasio@eio.uva.es}}
\and
\author{\fnms{Jean-Michel} \snm{Loubes}\ead[label=e2]{loubes@math.univ-toulouse.fr}}
%\and
%\author{\fnms{Third} \snm{Author}\thanksref{t1,m2}
%\ead[label=e3]{third@somewhere.com}
%\ead[label=u1,url]{http://www.foo.com}}

\thankstext{t1}{Research partially supported by Ministerio de Econom\'{\i}a y Competitividad, grant MTM2014-56235-C2-1-P}
%\thankstext{t2}{First supporter of the project}
%\thankstext{t3}{Second supporter of the project}
\runauthor{Del Barrio and Loubes}

\affiliation{IMUVA, Universidad de Valladolid and IMT, Universit\'e de Toulouse}

\address{Address of the First author\\
IMUVA, Universidad de Valladolid \\ Spain \\
\printead{e1}}

\address{Address of the Second author\\ IMT, Universit\'e de Toulouse\\France\\
\printead{e2}}
\end{aug}

\begin{abstract}
We consider the problem of optimal transportation with quadratic cost between a empirical 
measure and a general target probability on $\mathbb{R}^d$, with $d\geq 1$.
We provide new results on the uniqueness and stability of the associated optimal transportation
potentials, namely, the minimizers in the dual formulation of the optimal transportation problem.
As a consequence, we show that a CLT holds for the empirical transportation cost 
under mild moment and smoothness requirements. The limiting distributions are Gaussian and
admit a simple description in terms of the optimal transportation potentials.
\end{abstract}



\end{frontmatter}

\section{Introduction}

The analysis of the minimal transportation cost between two sets of random points or
of the transportation cost between an empirical and a reference measure is by now a classical
problem in probability, to which a significant amount of literature has been devoted. \\
\indent 
In the case of two sets of $n$ random points, say $X_1,\ldots,X_n$ and $Y_1,\ldots,Y_n$ in $\mathbb{R}^d$,
the object of interest is
$$T_{c,n}=\min_{\sigma}\frac 1 n\sum_{i=1}^n c(X_i,Y_{\sigma(i)}),$$
where $\sigma$ ranges is the set of permutations of $\{1,\ldots,n\}$ and $c(\cdot,\cdot)$ is some cost function. 
$T_{c,n}$ is usually referred to as the cost of optimal matching. This optimal matching problem is closely related
to the Kantorovich optimal transportation problem, which, in the Euclidean setting amounts to the minimization of
$$I[\pi]=\int_{\mathbb{R}^d\times \mathbb{R}^d} c(x,y) d\pi(x,y),$$
with $\pi$ ranging in the set of joint probabilities on $\mathbb{R}^d\times \mathbb{R}^d$
with marginals $P,Q$. Here $P$ and $Q$ are two probability measures on $\mathbb{R}^d$ and
the minimal value of $I[\pi]$ is known as the optimal transportation cost between $P$ and $Q$.
The cost functions $c(x,y)=\|x-y\|^p$ have received special attention and we will write
$\mathcal{W}_p^p(P,Q)$ for the optimal transportation cost in that case. It is well known that
with this choice of cost function $T_{c,n}=\mathcal{W}_p^p(P_n,Q_n)$, with $P_n$ and $Q_n$ denoting
the empirical measures on $P_n$ and $Q_n$. A related functional of interest is $\mathcal{W}_p^p(P_n,Q)$,
the transportation cost between the empirical measure on the sample $X_1,\ldots,X_n$ and a given probability
$Q$.

How large is the cost of optimal matching, $\mathcal{W}_p^p(P_n,Q_n)$? Under the assumption that $X_1,\ldots,X_n$
are i.i.d. with distribution $P$, $Y_1,\ldots,Y_n$ are i.i.d. with distribution  $Q$ and $P$ and $Q$ have finite $p$-th moment is is easy to conclude that
$\mathcal{W}_p^p(P_n,Q_n)\to \mathcal{W}_p^p(P,Q)$ almost surely. One might then wonder about the rate of approximation,
that is, how far is the empirical transportation cost from its theoretical counterpart. Much effort has been devoted
to the case when $P=Q$, namely, when the two random samples come from the same random generator. In this case 
$\mathcal{W}_p^p(P,Q)=0$ and the goal is to determine how fast does the empirical optimal matching cost vanish.
From the early work \cite{AjtaiKomlosTusnady}, followed by the important contributions \cite{Talagrand1992}, \cite{Talagrand}, 
\cite{TalagrandYukich} and \cite{DobricYukich}, it is known that the answer depends on the dimension 
$d$. In the case when $P=Q$ is the uniform distribution on the unit hypercube $\mathcal{W}_p(P_n,Q_n)=O(n^{-1/d})$,
if $d\geq 3$, with a slightly worse rate if $d=2$. The results for $d\geq 3$ were later extended to a more general setup
covering the case when $P=Q$ has bounded support and a density satisfying some smoothness
requirements. The one-dimensional case is different. If $p=1$ then, under some integrability assumptions
$\mathcal{W}_1(P_n,P)=O_P(n^{-1/2})$, with $\sqrt{n}\mathcal{W}_1(P_n,P)$ converging weakly to a non Gaussian limit,  see \cite{dBGM99}. 
If $p>1$ then it is still possible to get a limiting distribution for $\sqrt{n}\mathcal{W}_p(P_n,P)$, but now 
integrability assumptions are not enough and the available results require some smoothness conditions
on $P$ (and on its density), see \cite{dBGU05} for the case $p=2$. In fact, see \cite{BovkovLedoux2014}, 
the condition that $P$ has a positive density in 
an interval is necessary for boundedness of the sequence $\sqrt{n}E(\mathcal{W}_p(P_n,P))$ if $p>1$. In a different setting using PDE, rates in dimension 2 are also given~\cite{ambrosio2016pde}.

This paper provides CLT's and variance bounds for the quadratic transportation cost between an empirical measure
based on i.i.d. observations and a probability on $\mathbb{R}^d$ or between two sets of $d$-dimensional i.i.d. observations. 
More precisely, we will consider i.i.d. $\mathbb{R}^d$ valued random variables (r.v.'s in the sequel) 
$X_1,\ldots,X_n$ with common distribution $P$
and an additional probability $Q$ on $\mathbb{R}^d$. We will write $P_n$ for the empirical measure on $X_1,\ldots,X_n$
and will give CLT's for $\mathcal{W}_2(P_n,Q)$ (see our Theorem \ref{CLTdimd}). 
We also extend this result to CLT's for $\mathcal{W}_2(P_n,Q_m)$ when $Q_m$ is
the empirical measure on a further independent sample of i.i.d. r.v.'s, $Y_1,\ldots,Y_m$, with law $Q$.

Beyond the theoretical interest of the problem, we would like to emphasize the potential impact on 
statistical applications of our results. 
Quoting from \cite{SommerfeldMunk}, the transportation cost distance
`is an attractive tool for data analysis but statistical inference is hindered by the lack of distributional limits'.
This has led to some attempts to provide some distributional limits in different setups. In \cite{RipplMunkSturm} a
related (but different) problem is considered. There the sample $X_1,\ldots,X_n$ consists of i.i.d. Gaussian r.v.'s
(this is extended to cover elliptical models as well) and CLT's are given for the transportation between the
underlying Gaussian law and a Gaussian law with estimated parameters (see Theorems 2.1 and 2.2 there).
To our best knowledge the only work that deals with the issue of distributional limit laws for the transportation
cost between empirical measures is \cite{SommerfeldMunk} (see Theorem 1 there). However, the problem
considered there is of a different nature. Both generating probabilities $P$ and $Q$ are assumed to have finite support.
This allows to deal with the transportation cost as a functional of the multinomial vector of empirical frequencies, 
and the result follows from the directional Hadamard differentiability of this functional.
On the other hand we focus on the case when the probabilities $P$ and $Q$ are smooth (or at least one of them) 
and this requires the exploration of alternative methods of proof.

Our approach to the transportation cost between empirical measures comes from a closer
analysis of the Kantorovich duality. We give a self-contained description of this in Section \ref{Section2}
below. For the moment we limit ourselves to note that the transportation cost $\mathcal{W}_2^2(P,Q)$
can be expressed as 
$$\mathcal{W}_2^2(P,Q)=\int_{\mathbb{R}^d} \|x\|^2 dP(x)+\int_{\mathbb{R}^d} \|y\|^2 dQ(y)+2 \min_{(\varphi, \psi)\in \Phi} 
J(\varphi,\psi),$$
where $\Phi$ denotes the set of pairs of functions $(\varphi, \psi)\in L_1(P)\times L_1(Q)$ such that
$\varphi(x)+\psi(y)\geq x\cdot y$ and $J$ is the linear functional
$$J(\varphi,\psi)=\int_{\mathbb{R}^d}\varphi dP+\int_{\mathbb{R}^d}\psi dQ.$$
If $(\varphi,\psi)$ is a minimizing pair in $\Phi$ for $J$ we will refer
to $\psi$ as an \textit{optimal transportation potential} for the transportation
of $Q$ to $P$. The motivation for the name comes from the fact that, provided
$Q$ has a density, the optimal transportation problem is equivalent to the Monge transportation
problem, that is, 
$$\mathcal{W}_2^2(P,Q)=\min_{T:\, Q\circ T^{-1}=P}\int_{\mathbb{R}^d} \|x-T(x)\|^2 dQ(x),$$
(here and in the sequel $Q\circ T^{-1}$ denotes the law induced from $Q$ by the measurable map $T$).
A minimizing $T$ in the Monge problem is called an optimal transportation map from $Q$ to $P$. It 
is well known (see, e.g., Theorem~ 2.12 in \cite{Villani}) that the optimal transportation map is unique
if $Q$ has a density and, in fact, it is the unique map of the form $\nabla \psi$ with $\psi$
a proper, lower semicontinuous, convex function, that maps $Q$ to $P$ (see details below). It is also 
true that $\nabla \psi$ is an optimal transportation map if and only if $\psi$ is
an optimal transportation potential. Beyond uniqueness, it is also known that
optimal transportation maps enjoy some stability: if $P_n\to P$ in $\mathcal{W}_2^2$ distance
then the optimal transportation map from $Q$ to $P_n$ converges $Q$- almost surely to the optimal
transportation map from $Q$ to $P$, see for instance Corollary 5.23 in \cite{VillaniBig}. Optimal transportation potentials do not enjoy uniqueness or stability
in general. However, we show in this paper that they are essentially unique (up to the addition of a constant)
and that suitable versions can be chosen for which stability does hold. 

Once we have proved these stability results, our approach to the CLTs for the empirical transportation cost
relies on a rather simple application of the Efron-Stein variance inequality (see Section \ref{CLTs}).
Related techniques had been used to provide exponential concentration bounds for the empirical 
transportation cost (see \cite{delBarrioMatran2013a}). Here we give only variance bounds, but get two main
advantages. First, these variance bounds hold in great generality, requiring only finite fourth moments 
(in \cite{delBarrioMatran2013a} a bounded support is assumed for the exponential bounds). Second, they can be
adapted to prove a linearization result that yields as a direct consequence our CLT's that are presented in Section \ref{CLTsequel}.
The Efron-Stein method for variance inequalities boils down to bounding the
moments of the increase that results in replacing a member of a sample by an independent
copy. This is particularly convenient in optimal transportation, where a solution which is optimal for a sample
results in or can be transformed into a different solution which is not
optimal for the transformed sample, but yields a workable bound for the increase in transportation cost. We would
like to mention that we use this observation both for the primal and the dual formulation of the
transportation problem and that both uses are needed to prove our linearization result.

Finally, to end this introduction, we would like to explain the particularities that made ourselves
constrain our approach to the quadratic cost. While it was this quadratic case that historically 
received first a closer attention, the theory has then broadened and much of the key results have been
extended to more general costs. Of course, the Kantorovich duality holds in much greater generality.
Equivalence to the Monge version of optimal transportation requires, however, some additional
assumptions, related to strict convexity of the cost function. It does hold for the cost
$\|x-y\|^p$ with $p> 1$ and there are uniqueness and stability results for the optimal
transportation maps in this more general setup (see \cite{GangboMcCann}). However, our
approach to prove uniquess and stability of optimal transportation potentials
relies on tools from the theory of graphical convergence of multivalued maps (a particular
case of set convergence in the Painlev\'e-Kuratowski sense, see details in Section \ref{Section2}
below) which are particularly suited for the analysis of convex functions and their subgradients.
We expect that similar results will be developed to enable to handle version related to generalized 
concavity, which would allow to extend the approach in this paper to costs $\|x-y\|^p$ with $p> 1$.
This will be covered in a future work.

\section{Uniqueness and stability of optimal transportation potentials.}\label{Section2}

An essential component in our approach is the Kantorovich duality, which we succintly describe next and refer 
to the excellent monographs \cite{RachevRuschendorf}, \cite{Villani} or \cite{VillaniBig} for further details. 
Given Borel probabilities $P$ and $Q$ on $\mathbb{R}^d$ with finite second moment,
the optimal transportation problem (with quadratic cost) is the problem of minimization of 
$I[\pi]=\int_{\mathbb{R}^d\times\mathbb{R}^d } \|x-y\|^ 2 d\pi(x,y)$ in 
$\pi\in \Pi(P,Q)$, the set of Borel probability measures on $\mathbb{R}^d\times \mathbb{R}^d$
with marginals $P$ and $Q$. 
It is convenient to consider the equivalent problem of maximization of 
$\tilde{I}[\pi]=\int_{\mathbb{R}^d\times\mathbb{R}^d } x\cdot y d\pi(x,y)$ (note that $I[\pi]=\int \|x\|^2 dP(x)+
\int \|y\|^2 dQ(y)-2\tilde{I}[\pi]$).
We denote by $\Phi$ the set of pairs of functions $(\varphi, \psi)\in L_1(P)\times L_1(Q)$ such that
$$\varphi(x)+\psi(y)\geq x\cdot y$$
for every $x$ and $y$. We write also
\begin{equation}\label{dualfunctional}
J(\varphi,\psi)=\int_{\mathbb{R}^d}\varphi dP+\int_{\mathbb{R}^d}\psi dQ.
\end{equation}
Then,
\begin{equation}\label{KantorDual}
\min_{(\varphi,\psi)\in \Phi} J(\varphi,\psi)=\max_{\pi\in \Pi(P,Q)} \tilde{I}[\pi].
\end{equation}
With this result, to which we will refer as the Kantorovich duality, we are summarizing a 
number of different facts. First, the functional $\tilde{I}[\pi]$ admits a maximizer
in $\Pi(P,Q)$; second, the functional $J(\varphi,\psi)$ admits a minimizer in $\Phi$; finally, the optimal values are equal 
(see for instance Theorems 1.3 and 2.9 in\cite{Villani}). Furthermore, the maximizing pair for $J$, $(\varphi,\psi)$, 
can be taken to be 
a pair of lower semicontinuous, proper convex conjugate functions, that is, $\varphi(x)={\psi}^*(x)$, where 
$$h^*(x)=\sup_{y\in \mathbb{R}^d} (x\cdot y-h(y))$$
denotes the convex conjugate of $h$ (note that $\int {\varphi} dP+\int {\psi}dQ\geq \int {\psi}^* dP+\int \tilde{\psi}dQ$
since ${\varphi}\geq {\psi}^*$ if $(\varphi,\psi)\in \Phi)$). This results in a more precise description of the maximizers of $\tilde{I}[\pi]$, as follows.

For any $\pi\in \Pi(P,Q)$ and any $(\psi^*,\psi)$ in $\Phi$ we clearly have
$$J(\psi^*,\psi)=\int_{\mathbb{R}^d\times\mathbb{R}^d } (\psi^*(x)+\psi(y)) d\pi(x,y)\geq 
\int_{\mathbb{R}^d\times\mathbb{R}^d } x\cdot y d\pi(x,y)=\tilde{I}[\pi].$$
The Kantorovich duality (\ref{KantorDual}) entails that $(\psi^*,\psi)$ is a minimizer of $J$ and $\pi$ is a maximizer of
$\tilde{I}$ if and only if 
$$\int_{\mathbb{R}^d\times\mathbb{R}^d } (\psi^*(x)+\psi(y)-x\cdot y) d\pi(x,y)=0,$$
that is, if and only if the nonnegative function $\psi^*(x)+\psi(y)-x\cdot y$ vanishes $\pi$-almost surely. The condition
$\psi^*(x)+\psi(y)-x\cdot y=0$ holds if and only if $x\in\partial \psi(y)$ (if and only if $y\in\partial \psi^*(x)$). Here 
$\partial\psi (y)$ denotes the subgradient of $\psi$ at $y$, that can be written as 
$$\partial\psi(y)=\{ z\in\mathbb{R}^d: \psi(y')-\psi(y)\geq z\cdot (y'-y) \mbox{ for all } y'\in\mathbb{R}^d\},$$
which is a nonempty set if $\psi$ is a proper convex function and $y$ belongs to the interior of its domain (see \cite{Rockafellar} 
for further details). If $\psi$ is differentiable at $y$ then $\partial\psi (y)=\{\nabla \psi(y)\}$, where $\nabla$ denotes the usual 
gradient. We note that convex functions are locally Lipschitz, hence, by Rademacher's Theorem (see, e.g., p. 81 in
\cite{GariepyEvans}) they are differentiable at almost every point in the interior of their domain.
These facts can be used to prove that if $Q$ 
does not give mass to sets of Hausdorff dimension $d-1$ (in particular if $Q$
is absolutely continuous with respect to $\ell_d$, the Lebesgue measure on $\mathbb{R}^d$), then
(see Theorem 2.12 in \cite{Villani}) $(\psi^*,\psi)$ is a minimizing pair for $J$ if and only if 
$Q\circ (\nabla \psi)^{-1}=P$ and then $\pi=Q\circ(\nabla \psi,Id)^{-1}$ maximizes $\tilde{I}$.
The map $T=\nabla \psi$ 
is  known as the \textit{optimal transportation map} from $Q$ to $P$ and is $Q$-a.s. unique:
if $\psi_1$ were a further convex function such that $Q\circ (\nabla \psi_1)^{-1}=P$ then $\nabla \psi=\nabla \psi_1$
$Q$-almost surely.

Unlike the optimal transportation map, the \textit{optimal transportation potential}, that is a convex, lower semicontinuous $\psi$ 
such that $(\psi^*,\psi)$ minimizes $J$ (equivalently, a convex, lower semicontinuous $\psi$ 
such that $Q\circ (\nabla \psi)^{-1}=P$), is not unique, since, obviously $J(\psi^*-C,\psi+C)= 
J(\psi^*,\psi)$ for every $C\in\mathbb{R}$. However, under some additional regularity on $Q$ we can ensure that this
is the only way to produce a different optimal transportation potential. 
Our next result would be trivial if we were imposing further smoothness assumptions on the convex potentials:
two differentiable functions on a convex domain that have the same gradient are equal up 
to addition of a constant. What we show next is that, for convex functions, having a common gradient at \textit{almost} every point is enough to reach
the same conclusion. 
\begin{Lemma}\label{uniquepotential}
Assume $\psi_1$ and $\psi_2$ are finite convex functions on a nonempty convex, open set $A\subset \mathbb{R}^d$ such that
$$\nabla \psi_1(x)=\nabla \psi_2(x)\quad \mbox{ for almost every } x\in A.$$
Then there exists $C\in \mathbb{R}$ such that $\psi_1(x)=\psi_2(x)+C$ for all $x\in A$.
\end{Lemma}

\medskip
\noindent
\textbf{Proof.} For $i=1,2$, we write $\partial \varphi_i(x)$ for the subgradient of $\varphi_i$ at $x\in A$, namely, the set of $z\in\mathbb{R}^d$ such
that $\varphi_i(y)-\varphi_i(x)\geq z\cdot (y-x)$ for all $y\in\mathbb{R}^d$. 
We also write $S_i(x)$ for the set of points $z\in\mathbb{R}^d$ such that $z=\lim_{n\to\infty} \nabla \varphi_i(x_n)$ for some sequence
$x_n$ which satisfies $\lim_{n\to\infty}x_i=x$. Then $\partial \varphi_i(x)$ is the closure of the convex hull of $S_i(x)$ (see Theorem 25.6, p. 246
in \cite{Rockafellar}; note that the normal cone to a point in the interior of the domain of a convex function is simply $\{0\}$).  
Now, assume that $z\in S_1(x)$, with $z=\lim_{n\to\infty} \nabla \varphi_1(x_n)$ and $x_n$ is some sequence converging to $x\in A$.
Denote by $B\subset A$ the set such that $A-B$ has null Lebesgue measure while for $x\in B$ $\varphi_i$, $i=1,2$ are differentiable 
at $x$ with $\nabla\varphi_1(x)=\nabla \varphi_2(x)$. We note that $\nabla \varphi_1$ is continuous in the set of points of 
differentiability of $\varphi_1$ (Theorem 25.5 in \cite{Rockafellar}). Hence, for each $n$ we can find $\tilde{x}_n\in B$ such 
that $\|x_n-\tilde{x}_n\|\leq \frac 1 n$ and $\|\nabla \varphi_1(x_n)-\nabla \varphi_1(\tilde{x}_n)\|\leq \frac 1 n$. But then
$\tilde{x}_n\to x$ and $\nabla \varphi_1(\tilde{x}_n)= \nabla \varphi_2(\tilde{x}_n)\to z$, which shows that $z\in S_2(x)$ and implies that
$S_1(x)\subset S_2(x)$. By symmetry, we also have $S_2(x)\subset S_1(x)$. Now, two convex functions with equal subgradient at every point
are equal up to the addition of a constant (see Theorem 24.9 in \cite{Rockafellar}; we note that although the statement of this Theorem considers
convex functions on $\mathbb{R}^d$ the proof can be reproduced verbatim for convex functions on a smaller 
convex, open domain in $\mathbb{R}^d$). This completes the proof.
\hfill $\Box$

\medskip
As a consequence of Lemma \ref{uniquepotential}, we obtain uniqueness of optimal transportation potentials (up to the addition of a constant) 
under suitable regularity assumptions.

\begin{Corollary}\label{uniqueotp}
Assume that $P$ and $Q$ are Borel probabilities on $\mathbb{R}^d$ with finite second moments and
\begin{equation}\label{regularQ}
Q \mbox{ has a positive density in the interior of its convex support.} 
\end{equation}
Then, if $\psi_1$, $\psi_2$ are convex, lower semicontinuous convex functions such that $J(\psi_1^*,\psi_1)=
J(\psi_2^*,\psi_2)=\min_{(\varphi,\psi)\in \Phi} J(\varphi,\psi)$, with $J(\varphi,\psi)$ as in 
(\ref{dualfunctional}), there exists $C\in \mathbb{R}$ such that $\psi_2=\psi_1+C$ in the interior of the support of $Q$.
In particular, $\psi_2=\psi_1+C$ $Q$-a.s..
\end{Corollary}

\medskip
\noindent
\textbf{Proof.} Uniqueness of the optimal transportation map and (\ref{regularQ}) ensure that $\nabla {\psi}_1(x)=\nabla {\psi}_2(x)$
for almost every $x\in A$, the interior of the support of $Q$. Lemma \ref{uniquepotential} allows to conclude that
$\psi_2(x)={\psi}_1(x)+C$ for some constant $C$ and every $x$ in the interior of $A$. The conclusion
follows from the fact that the boundary of a convex set has zero Lebesgue measure. 
%Uniqueness of the optimal transportation map and (\ref{regularQ}) ensure that $\nabla {\varphi}_1(x)=\nabla {\varphi}_2(x)$
%for almost every $x\in A$, the interior of the support of $Q$. Lemma \ref{uniquepotential} allows to conclude that
%$\varphi_2(x)={\varphi}_1(x)+C$ for some constant $C$ and every $x$ in the interior of $A$. The conclusion
%follows from the fact that the boundary of a convex set has zero Lebesgue measure. 

\hfill $\Box$

\begin{Remark}{\em
Uniqueness of the optimal transportation potential fails without assumption (\ref{regularQ}). As a counterexample,
consider the probability $P$ giving mass $\frac 1 2$ to the points $-1,1$ and assume that $Q_\varepsilon$ is the uniform law 
on the set $(-\varepsilon-1,-\varepsilon)\cup (\varepsilon,1+\varepsilon)$, $\varepsilon>0$. Non-decreasing 
maps are optimal. Hence, the optimal transportation map from $Q_\varepsilon$ to $P$ is $T_\varepsilon(x)=-1$, 
$x<0$, $T_\varepsilon(x)=1$, $x>0$. The maps $\psi_{\varepsilon,L}(x)=-x$, $x\leq -\frac L 2$,  
$\psi_{\varepsilon,L}(x)=x+L$, $x\geq -\frac L 2$, $0<L<\varepsilon$, are continuous, convex and
satisfy $\psi_{\varepsilon,L}'=T_\varepsilon$ $Q_\varepsilon$ a.s. . Hence, they are optimal transportation
potentials. However, if $L_1\ne L_2$, then there is no choice of a constant $C$ such that $\psi_{\varepsilon,L_2}=\psi_{\varepsilon,L_1}+C$
$Q_\varepsilon$ a.s. This example can be easily adapted to general dimension.

\hfill $\Box$
}
\end{Remark}

We turn now to stability in optimal transportation problems. 
We will assume that $Q$ is a regular probability measure on $\mathbb{R}^d$ (in the sense of
(\ref{regularQ})) and $P_n,P$ are probabilities satisfying $\mathcal{W}_2(P_n,P)\to 0$. It is well known 
(see Theorem 3.4 in \cite{CuestaAlbertosMatranTueroDiaz}) that the optimal transportation maps from $Q$ to $P_n$, say $T_n$, converge $Q$-a.s. to $T$,
the optimal transportation map from $Q$ to $P$. Here we will provide stability results for the
optimal transportation potentials.

 A main tool in our approach will be the concept of \textit{graphical convergence}
of multivalued maps, which is a particular case of set convergence in the Painlev\'e-Kuratowski sense. 
We include next a brief summary of some related key facts and refer to \cite{RockafellarWets} for a detailed
account of the main results on the topic.

Given a sequence of subsets $\{C_n\}_{n\geq 0}$ of $\mathbb{R}^d$, its outer limit, to be denoted
$\limsup_{n\to\infty} C_n$ is the set of points $x\in \mathbb{R}^d$ such that $x=\lim_{j\to\infty} x_{n_j}$
for some subsequence $n_j$ and some choice of points $x_{n_j}\in C_{n_j}$, while the inner
limit (denoted $\liminf_{n\to\infty} C_n$) is the set of points $x\in \mathbb{R}^d$ such that $x=\lim_{n\to\infty} x_{n}$
for some sequence $x_n$ such that $x \in C_n$ for all $n\geq n_0$ (for some $n_0$). Obviously,
$\liminf_{n\to\infty} C_n\subset \limsup_{n\to\infty} C_n$. When these two sets are equal (to $C$, say)
then the sequence $C_n$ is said to converge to $C$ in the Painlev\'e-Kuratowski sense. The limiting sets are necessarily closed
and, in fact, it makes no difference to replace $C_n$ by its closure in all these definitions (see Proposition 
4.4 in \cite{RockafellarWets}).

A multivalued map, $T$, from $\mathbb{R}^d$ to $\mathbb{R}^d$ is a map that assigns to each $x\in\mathbb{R}^d$, a set
$T(x)\subset \mathbb{R}^d$. The domain of $T$ is the set of $x\in\mathbb{R}^d$ such that $T(x)\ne \emptyset$, while the graph
is the subset 
$$\mbox{gph}(T)=\Big\{(x,t)\in\mathbb{R}^d\times \mathbb{R}^d:\, t\in T(x)\Big\}.$$
Multivalued maps can be identified with subsets of $\mathbb{R}^d\times \mathbb{R}^d$. Given a set $T\subset
\mathbb{R}^d\times \mathbb{R}^d$ we can define the map $\tilde{T}$ by the 
rule $\tilde{T}(x)=\{t\in\mathbb{R}^d:\, (x,t)\in T\}$ and then the graph of $\tilde{T}$ equals $T$. This identification
allows to define convergence of multivalued maps in terms of set convergence of their
graphs in the Painlev\'e-Kuratowski sense. More precisely, the sequence of multivalued maps $\{T_n\}_{n\geq 1}$ from
$\mathbb{R}^d$ to $\mathbb{R}^d$ is said to converge graphically to $T$ if the graphs $\mbox{gph}(T_n)$ converge
to $\mbox{gph}(T)$ in the Painlev\'e-Kuratowski sense, see Chapter 5 in \cite{RockafellarWets} for details.
For convenience, we include next two results about convergence of sets and multivalued maps. The first one is 
a characterization of graphical convergence, which is just a rewriting
of Proposition 5.33 in \cite{RockafellarWets}. The second is a key result on sequential compactness in the 
Painlev\'e-Kuratowski sense. 

\begin{Proposition}\label{graphicalconvergence}
The sequence of multivalued maps $\{T_n\}_{n\geq 1}$ converges graphically to $T$ if and only if for every
$x\in\mathbb{R}^d$ the following two conditions hold:
\begin{itemize}
\item[(a)] if $x_n\to x$, $y_n\in T_n(x_n)$ for large $n$ and there is a subsequence $y_{n_j}\to y$, then $y\in T(x)$,
\item[(b)] if $y\in T(x)$ then there exist sequences $\{x_n\}$, $\{y_n\}$ with $x_n\to x$, $y_n\in T_n(x_n)$ for 
large $n$ and such that $y_n\to y$.
\end{itemize}
\end{Proposition}

\begin{Theorem}\label{compactness}

\begin{itemize}
\item[(a)] Assume that $\{C_n\}_{n\geq 1}\subset \mathbb{R}^d$ satisfies that for some $\varepsilon>0$ and some subsequence $\{n_j\}$  
$C_{n_j}\cap B(0,\varepsilon)\ne \emptyset$ for every $j\geq 1$, where $B(0,\varepsilon)$ denotes the open ball of 
radius $\varepsilon$ centered at the origin. Then there exists a subsequence $\{n_{j_k}\}$ and a nonempty subset $C\subset \mathbb{R}^d$ 
such that $C_{n_{j_k}}$ converges to $C$ in the Painlev\'e-Kuratowski sense.
\item[(b)] Assume that $\{T_n\}_{n\geq 1}$ is a sequence of multivalued maps from $\mathbb{R}^d$ to $\mathbb{R}^d$ such that for some 
bounded sets $C,D\subset \mathbb{R}^d$ and some subsequence $\{n_j\}$ there exist $x_{n_j}\in C$ with $T_{n_j}(x_{n_j})\cap D\ne \emptyset$
for all $j\geq 1$. Then there exists a subsequence $\{n_{j_k}\}$ and a multivalued map, $T$, from $\mathbb{R}^d$ to $\mathbb{R}^d$, 
with nonempty domain such that $T_{n_{j_k}}$ converges graphically to $T$.
\end{itemize}

\end{Theorem}

\medskip
\noindent \textbf{Proof.} We note that the assumption in (a) is simply a rewriting of the assumption in Theorem 4.18 in
\cite{RockafellarWets} (the condition that the sequence of sets does not \textit{escape to the horizon}). Similarly,
(b) follows from Theorem 5.36 in \cite{RockafellarWets}.

\hfill $\Box$

\medskip
The link between optimal transportation and the theory of multivalued maps comes from the fact that 
a transportation plan $\pi$ is optimal (a minimizer for $\tilde{I}$) if and only its support is contained
in the graph of the multivalued map $\partial \psi$ for some proper, lower semicontinuous, convex $\psi$ (recall
the discussion above; see also Theorem 2.12 in \cite{Villani}). It is well known that subgradients 
of convex maps can be characterized in terms of monotonicity or cyclical monotonicity. A multivalued map
$T$ from $\mathbb{R}^d$ to $\mathbb{R}^d$ is monotone if $(t_1-t_0)\cdot (x_1-x_0)\geq 0$ whenever
$t_i\in T(x_i)$, $i=0,1$. It is cyclically monotone if for every choice of $m\geq 1$, points $x_0,\ldots,x_m$
and elements $t_i\in T(x_i)$, $i=0,\ldots,m$, we have
$$t_0\cdot (x_1-x_0)+t_1\cdot (x_2-x_1)+\cdots+t_m\cdot (x_0-x_m)\leq 0.$$
A monotone multivalued map is maximal monotone if its graph cannot be enlarged without losing the 
monotonicity property and similarly for maximal cyclically monotone maps. It is easy to see that 
every cyclically monotone map is also monotone. It is also true that a maximal cyclically monotone map is
maximal monotone and, in fact, a multivalued map $T$ has the form $T=\partial \psi$ for some proper, lower 
semicontinuous, convex $\psi$ if and only if $T$ is maximal cyclically monotone 
(see Theorems 12.17 and 12.25 in \cite{RockafellarWets}; Theorem 12.25 is often referred to as `Rockafellar's Theorem').

\medskip
In our stability result for optimal transportation potential we will make use of the following result on convergence
of cyclically monotone maps. While it follows easily from related known results, we have not been able to find
it in the literature and therefore states its result in the following theorem.

\begin{Theorem}\label{stablecyclic}
If a sequence of cyclically monotone maps $\{T_n\}$ from $\mathbb{R}^d$ to $\mathbb{R}^d$ converges graphically
then the limit map, $T$, must be cyclically monotone. If the $T_n$ are maximal cyclically monotone then
$T$ is also maximal cyclically monotone.

Assume $\{\psi_n\}$ is a sequence of proper, lower semicontinuous, convex maps from $\mathbb{R}^d$ to $\mathbb{R}$ such that
for some bounded sets $C,D\subset \mathbb{R}^d$ and some subsequence $\{n_j\}$ there exist $x_{n_j}\in C$ with $\partial \psi_{n_j}(x_{n_j})\cap D\ne \emptyset$
for all $j\geq 1$. Then there exists a subsequence $\{n_{j_k}\}$ and a proper, lower semicontinuous, convex map, $\psi$, 
from $\mathbb{R}^d$ to $\mathbb{R}$, with subgradient with nonempty domain such that $\partial \psi_{n_{j_k}}$ converges graphically to $\partial \psi$.
\end{Theorem}

\medskip
\noindent \textbf{Proof.} Take $t_i\in T(x_i)$, $i=0,\ldots,m$. The points $(x_i,t_i)$ belong to the graph of $T$, hence they 
belong to $\liminf_{n\to \infty}\mbox{gph}(T_n)$ and, consequently, there are sequences $(x_{n,i},t_{n,i})\in \mbox{gph}(T_n)$ (for large enough
$n$) such that $(x_{n,i},t_{n,i})\to (x_i,t_i)$, $i=0,1,\ldots,m$. By cyclical monotonicity we have 
$$t_{n,0}\cdot (x_{n,1}-x_{n,0})+t_{n,1}\cdot (x_{n,2}-x_{n,1})+\cdots+t_{n,m}\cdot (x_{n,0}-x_{n,m})\leq 0.$$
Taking limits we conclude that 
$$t_0\cdot (x_1-x_0)+t_1\cdot (x_2-x_1)+\cdots+t_m\cdot (x_0-x_m)\leq 0.$$
Therefore $T$ is cyclically monotone. If $T_n$ are maximal cyclically monotone then they are maximal monotone. By Theorem
12.32 in \cite{RockafellarWets} $T$ must be maximal monotone. Hence, it is also maximal cyclically monotone 
(if we could enlarge the graph of $T$ preserving cyclical monotonicity, then the enlarged graph would also be monotone,
contradicting maximal monotonicity).

For the second part we use Rockafellar's theorem and part (b) of Theorem \ref{compactness}.

\hfill $\Box$

Finally, we quote a technical result relating graphical convergence of subgradients of convex functions to
pointwise convergence of the convex functions themselves. A proof follows easily from Theorem 12.35 and
Exercise 12.36 in \cite{RockafellarWets}.

\begin{Proposition}\label{convexpotentialconvergence}
Assume $\psi$, $\{\psi_n\}$ are proper, lower semicontinuous, convex maps from $\mathbb{R}^d$ to $\mathbb{R}$ such that
$\partial \psi_n$ converges to $\partial \psi$ graphically and there is a sequence $(x_n,t_n)$ with $t_n\in\partial \psi_n(x_n)$
and a pair $(x_0,t_0)$ with $t_0\in\partial \psi(x_0)$ satisfying $(x_n,t_n)\to (x_0,t_0)$ and $\psi_n(x_n)\to \psi (x_0)$. Then, if $\psi$ is finite at $x$, 
$\tilde{x}_n\to x$ and $\liminf_{n\to \infty} \partial \psi_n(\tilde{x}_n)\ne \emptyset$ we have
$$\lim_{n\to\infty}\psi_n(\tilde{x}_n)=\psi(x).$$
\end{Proposition}

We are now ready for the announced result on stability of optimal transportation potentials.

\begin{Theorem}\label{convergenceofpotentials}
Assume $Q$ satisfies (\ref{regularQ}) and $Q_n$, $P_n$ ,$P$ are probabilities such that $\mathcal{W}_2(P_n,P)\to 0$ and $\mathcal{W}_2(Q_n,Q)\to 0$. 
If $\psi_n$ (resp. $\psi$) are optimal transportation potentials from $Q_n$ to $P_n$ (resp. from $Q$ to $P$) then there exist constants $a_n$ such that if 
$\tilde{\psi}_n=\psi_n-a_n$ then $\tilde{\psi}_n(x)\to \psi(x)$ for every $x$ in the interior of the support of $Q$, hence, for $Q$-almost every $x$.
\end{Theorem}

\medskip
\noindent \textbf{Proof.} We write $\pi_n$ for an optimal transportation plan for $Q_n,P_n$ and $\pi$ for the optimal transportation plan for $Q,P$.
We recall that $\pi$ is unique and $\pi=Q\circ (Id, \nabla \psi)^{-1}$. $\pi$ is concentrated in the graph of $\partial \psi$, that is, in the closed set
$\{(x,y)\in\mathbb{R}^d\times \mathbb{R}^d: \psi(x)+\psi^*(y)=x\cdot y\}=
\{(x,y)\in\mathbb{R}^d\times \mathbb{R}^d: y\in \partial\psi(x)\}$. It is easy to see that $\pi_n\to \pi$ weakly. As before, we denote by
$A$ the interior of the support of $Q$. 
We write $\tilde{A}$ for the set of $x\in A$ such that $\psi$ is differentiable at $a$. 
Then $Q(\tilde{A})=Q(A)=1$. Furthermore (see Theorem 25.5 in \cite{Rockafellar})
$\nabla \psi$ is continuous at every differentiability point $x\in A$. Fix $x_0\in \tilde{A}$ and set $y_0=\nabla\psi(x_0)$.
Now, for every $\varepsilon>0$ there exists $\delta>0$ such that $\|\nabla\psi(x)-y_0\|\leq \varepsilon$ if $x\in \tilde{A}$
and $\|x-x_0\|\leq \delta$. Hence, $\pi(B(x_0,\delta)\times B(y_0,\varepsilon))\geq Q(B(x_0,\delta))=\eta >0$ by Assumption
(\ref{regularQ}), and weak convergence implies that $\pi_n (B(x_0,\delta)\times B(y_0,\varepsilon))\geq \frac \eta 2$ for large enough
$n$. But $\pi_n$ is concentrated in the graph of $\partial \psi_n$, hence, there exists $(x_n,y_n)$ with $y_n\in\partial \psi_n(x_n)$,
$\|x_n-x_0\|<\delta$, $\|y_n-y_0 \|<\varepsilon$. 
We take now a sequence $\varepsilon_k\searrow 0$. For every $k\geq 1$ we choose
$\delta_k\in(0,\frac 1 k)$ such that $\|\nabla \psi (x)-y_0\|<\varepsilon_k$ if $\|x-x_0\|<\delta_k$ and $x\in \tilde{A}$.
As before, $\pi(B(x_0,\delta_k)\times B(y_0,\varepsilon_k))\geq Q(B(x_0,\delta_k))=\eta_k>0$. Fix $n_0=0$ and, for $k\geq 1$,
$n_k>n_{k-1}$ such that $\pi_n(B(x_0,\delta_k)\times B(y_0,\varepsilon_k))\geq \frac{\eta_k}2$ if $n\geq n_k$. Recall that 
$\pi_n$ is concentrated in the graph of $\partial \psi_n$. For $n=1,\ldots,n_1-1$ we take any pair $(x_n,y_n)$ with $y_n\in\partial 
\psi_n(x_n)$. For $k\geq 2$ and $n=n_{k-1},\ldots,n_k$ we take $(x_n,y_n)\in B(x_0,\delta_{k-1})\times B(y_0,\varepsilon_{k-1})$
such that $y_n\in\partial \psi_n(x_n)$. This construction yields a sequence $(x_n,y_n)$ such that $x_n\to x_0$, $y_n\to y_0$
and $y_n\in\partial \psi_n(x_n)$. 
We note that $(x_0,y_0)\in\limsup \mbox{gph } \partial \psi_n$.
We set now $a_n=\psi_n(x_n)-\psi(x_0)$ and define $\tilde{\psi}_n(x)=\psi_n(x)-a_n$. Obviously,
$\partial \tilde{\psi}_n(x)=\partial {\psi}_n(x)$ for every $x$. By Theorem \ref{stablecyclic} there exists a proper, lower semicontinuous
convex function $\rho$ such that $\partial \tilde{\psi}_n$ converges graphically to $\partial \rho$ along a subsequence. We keep 
the same notation for the subsequence. We see that $y_0\in\partial \rho(x_0)$. We can consider now $x\in\mathcal{A}$, $y=\nabla \psi(x)$ and apply
the same argument to conclude that $y\in\partial\rho (x)$. This implies that $\mbox{dom }(\rho)\supset A$. Hence $\rho$ must be differentiable 
and $\nabla \rho(x)=\nabla \psi(x)$ at almost every point in $A$. We conclude, using Lemma \ref{uniquepotential}, that $\rho=\psi+C$ in $A$, hence,
subtracting a constant, if necessary, $\rho=\psi$ in $A$. Since $\tilde{\psi}_n(x_n)=\psi(x_0)=\rho(x_0)$, applying Proposition \ref{convexpotentialconvergence}
we obtain that $\tilde{\psi}_n(x)\to\rho(x)=\psi(x)$ for all $x\in\tilde{A}$, hence (see Theorem 7.17 in \cite{RockafellarWets}) 
$\tilde{\psi}_n(x)\to\rho(x)=\psi(x)$ for all $x\in A$. Note that from this argument we see, in fact, that for any $x\in A$ and any subsequence $n'$ we
can extract a further subsequence $n''$ such that $\tilde{\psi}_{n''}\to \psi(x)$. But this proves that $\tilde{\psi}_{n}\to \psi(x)$
as $n\to\infty$ for every $x\in A$. This completes the proof.

\hfill $\Box$

\begin{Remark} \label{generalization}{\em 
Theorem \ref{convergenceofpotentials} extends known results about stability of optimal transportation maps. 
In fact, it covers the case $Q_n=Q$. In this case $\psi_n$ is differentiable at almost every $x\in A$. From
the proof of Theorem \ref{convergenceofpotentials} we have graphical convergence of $\partial \psi_n$ to $\partial \rho$
with $\rho=\psi$ in $A$. This implies (see, e.g., Exercise 12.40 (a) in \cite{RockafellarWets}) that
$\nabla \psi_n(x)\to \nabla \psi(x)$ at almost every $x\in A$, that is $\nabla \psi_n\to \nabla \psi$ $Q$-a.s.. This stability
result for optimal transportation maps is contained in Theorem 3.4 in \cite{CuestaAlbertosMatranTueroDiaz} or in \cite{Heinich}. Our result applies
to a non-smooth setup in that the $Q_n$'s are not assumed to have a density (on the other hand, we need to impose
additional regularity assumptions on $Q$ to ensure convergence of the convex potentials).
}
\end{Remark}

\hfill $\Box$

Under some moment assumptions the stability result in Theorem \ref{convergenceofpotentials} can be complemented with 
$L_2$ convergence. As in the Introduction, in our next result 
$\mathcal{W}_4$ denotes the transportation cost metric associated to the cost function $c(x,y)=\|x-y\|^4$.
We note that the condition $\mathcal{W}_4(P_n,P)\to 0$ implies the weaker assumption $\mathcal{W}_2(P_n,P)\to 0$
and also that the conclusions in Theorem \ref{secondmomentth} do not depend on the particular choice of the potential $\psi$
since all the possible choices are $Q$-a.s. equal up 
to the addition of a constant.

\begin{Theorem}\label{secondmomentth}
Assume that $Q,P,\{P_n\}_{n\geq 1}$ are probabilities on $\mathbb{R}^d$ with finite fourth moment with $Q$
satisfying (\ref{regularQ}) and write $\psi$ (resp. $\psi_n$) for a proper, lower semicontinuous function such that
$\nabla \psi$ (resp. $\nabla \psi_n$) is the optimal transportation map from $Q$ to $P$ (resp. from $Q$ to $P_n$).
Then $\psi, \psi_n\in L_2(Q)$. Furthermore, if $\mathcal{W}_4(P_n,P)\to 0$, then taking $\tilde{\psi}_n$ as in Theorem
\ref{convergenceofpotentials} we have that $\tilde{\psi}_n\to \psi$ in $L_2(Q)$.
\end{Theorem}

\medskip
\noindent \textbf{Proof.} We keep the notation for $A$ and $\tilde{A}$ as in the proof of Theorem \ref{convergenceofpotentials} and
the choice of $x_0\in U$ and write $z_0=\nabla \psi(x_0)$.
Then 
\begin{equation}\label{lower}
\psi(x)\geq \psi(x_0)+z_0\cdot (x-x_0), \quad x\in\mathbb{R}^d.
\end{equation}
On the other hand, since $z_0\in\partial \psi(x_0)$ we have $\psi(x_0)+\psi^*(z_0)=x_0\cdot z_0$, hence, 
$x_0\in\partial \psi^*(z_0)$ and .  
\begin{equation}\label{lower2}
\psi^*(z)\geq \psi^*(z_0)+x_0\cdot (z-z_0), \quad z\in\mathbb{R}^d.
\end{equation}
But optimality implies that $\psi(x)+\psi^*(\nabla \psi(x))=x \cdot \nabla \psi (x)$ $Q$-a.s.. Therefore,
using (\ref{lower2}) we conclude that, $Q$-a.s.,
\begin{equation}\label{upper}
\psi(x)\leq x\cdot \nabla \psi(x) - \psi^*(z_0)-x_0\cdot (\nabla \psi(x)-z_0)=\psi(x_0)+(x-x_0)\cdot \nabla \psi(x).
\end{equation}
Combining (\ref{lower}) and (\ref{upper}) we see that 
\begin{eqnarray*}
|\psi(x)-\psi(x_0)|&\leq & |(x-x_0)\cdot \nabla \psi (x_0)|+|(x-x_0)\cdot \nabla \psi (x)|\\
&\leq & \|x-x_0\|^2+{\textstyle  \frac 1 2} \|\nabla \psi (x_0)\|^2+{\textstyle  \frac 1 2} \|\nabla \psi (x)\|^2, \quad Q-\mbox{a.s.}
\end{eqnarray*}
By assumption $\|x-x_0\|^2$ is in $L_2(Q)$. Also, since, $\nabla \psi$ transports $Q$ to $P$, $\int \|\nabla \psi(x)\|^4 dQ(x)=\int \|z\|^4 dP(z)$.
This shows that $\psi\in L_2(Q)$. The same argument works for $\psi_n$ or $\tilde{\psi}_n$, in fact, 
\begin{eqnarray*}
|\tilde{\psi}_n(x)-\tilde{\psi}_n(x_0)|\leq  \|x-x_0\|^2+{\textstyle  \frac 1 2} \|\nabla \psi_n (x_0)\|^2+{\textstyle  \frac 1 2} \|\nabla \psi_n (x)\|^2, \quad Q-\mbox{a.s.}
\end{eqnarray*}
Now, $\|\nabla \psi_n (x)\|^4\to \|\nabla \psi(x)\|^4$ $Q$-a.s. and $\int \|\nabla \psi_n (x)\|^4dQ(x) \to \int \|\nabla \psi(x)\|^4dQ(x)$. Hence,
the sequence $\|\nabla \psi_n \|^4$ converges to $\|\nabla \psi \|^4$ in $L^1(Q)$ according to Scheff\'e Lemma. So it is $Q$-uniformly integrable, and the same applies to $\tilde{\psi}_n^2$, which combined with Theorem
\ref{convergenceofpotentials} proves that $\tilde{\psi}_n\to \psi$ in $L_2(Q)$. 

\hfill $\Box$

\section{Variance bounds.}\label{CLTs}
We turn now to concentration bounds and Central Limit Theorems for the empirical $L_2$-Wasserstein distance on $d$-dimensional data. 
From this point we assume that $P_n$ denotes the empirical measure on $X_1,\ldots,X_n$, i.i.d. r.v.'s with distribution $P$ and 
$P$ and $Q$ are Borel probabilities on $\mathbb{R}^d$ with finite second moments. 
A main tool in our proofs is the Efron-Stein inequality for variances, namely, that if $Z=f(X_1,\ldots,X_n)$ with
$X_1,\ldots,X_n$ independent random variables, $(X_1',\ldots,X_n')$ is an independent copy of $(X_1,\ldots,X_n)$ and
$Z_i=f(X_1,\ldots,X_i',\ldots,X_n)$ then
$$\mbox{Var}(Z)\leq \frac 1 2 \sum_{i=1}^nE(Z-Z_i)^2= \sum_{i=1}^n E(Z-Z_i)_+^2.$$
We refer, for instance, to \cite{BoucheronLugosiMassart2013} for a proof. In the particular case when $X_1,\ldots,X_n$ are
i.i.d. and $f$ is a symmetric function of $x_1,\ldots,x_n$ all the values $E(Z-Z_i)_+^2$ are equal and the bound simplifies
to 
\begin{equation}\label{EFsimple}
\mbox{Var}(Z)\leq n E(Z-Z')_+^2
\end{equation}
with $Z'=f(X_1',X_2,\ldots,X_n)$.

\medskip
We show first a variance bound for $\mathcal{W}_2^2(P_n,Q)$.

\begin{Theorem}\label{EfronSteindimd}
If $Q$ has a density and $P$ and $Q$ have finite fourth moments then
$$\mbox{\em Var}(\mathcal{W}_2^2(P_n,Q))\leq \frac{C(P,Q)}{n},$$
where $C(P,Q)=8\Big(E(\|X_1-X_2\|^2 \|X_1\|^2)+(E\|X_1-X_2\|^4)^{1/2} \Big({\textstyle \int_{\mathbb{R}^d}\|y\|^4dQ(y) } \Big)^{1/2} \Big)$.
\end{Theorem}

\medskip
\noindent
\textbf{Proof.} We write $Z=\mathcal{W}_2^2(P_n,Q)$. The assumption that $Q$ has a density ensures the existence
of an otimal transportation map, $T$, from $Q$ to $P_n$. Hence, denoting $C_i=\{y\in\mathbb{R}^d:\, T(y)=X_i\}$ we
have $Q(C_i)=\frac 1 n$ and 
$$Z=\sum_{i=1}^n \int_{C_i} \|y-X_i\|^2 dQ(y).$$
Let us consider an additional random variable $X_1'$ with law $P$, independent of $X_1,\ldots,X_n$, write $P_n'$ for the
empirical measure on $X_1',X_2,\ldots,X_n$ and $Z'=\mathcal{W}_2^2(P'_n,Q)$. Let us also denote by $T'$ the o.t.m.
from $Q$ to $P_n'$ and $C_1'=\{y\in\mathbb{R}^d:\, T'(y)=X_1'\}$, $C_i'=\{y\in\mathbb{R}^d:\, T'(y)=X_i\}$, $i=2,\ldots,n$.
Then 
$$Z'=\int_{C_1'} \|y-X_1'\|^2 dQ(y)+\sum_{i=2}^n \int_{C_i'} \|y-X_i\|^2 dQ(y),$$
while
$$Z\leq \int_{C_1'} \|y-X_1\|^2 dQ(y)+\sum_{i=2}^n \int_{C_i'} \|y-X_i\|^2 dQ(y).$$
This implies that
$$Z-Z'\leq \int_{C_1'} (\|y-X_1\|^2- \|y-X_1'\|^2)dQ(y)\leq \|X_1-X_1'\| \Big(\frac 1 n \big(\|X_1\|+\|X_1'\|\big)+2\int_{C_1'}\|y\|dQ(y)\Big),
$$
from which we conclude that
\begin{equation}\label{bound1}
E(Z-Z')_+^2\leq \frac 8 {n^2} E(\|X_1-X_1'\|^2 \|X_1\|^2)+ 8 E\Big(\|X_1-X_1'\|^2 \Big(\int_{C_1'}\|y\|dQ(y)\Big)^2 \Big).
\end{equation}
We note now that
$$\int_{C_1'}\|y\|dQ(y)\leq \Big(\int_{C_1'} 1 dQ(y)\Big)^{3/4}\Big(\int_{C_1'}\|y\|^4dQ(y)\Big)^{1/4}=
\frac 1  {n^{3/4}}\Big(\int_{C_1'}\|y\|^4dQ(y)\Big)^{1/4} .$$
By exchangeability we have $\int_{C_1'}\|y\|^4dQ(y)\overset d =\int_{C_1}\|y\|^4dQ(y)\overset d = \int_{C_j}\|y\|^4dQ(y)$,  for all $j=2,\ldots,n$.
This shows that
$$E\Big(\int_{C_1'}\|y\|^4dQ(y)\Big)=\frac 1 n E\Big(\sum_{j=1}^n\int_{C_j}\|y\|^4dQ(y)\Big)=\frac 1 n \int_{\mathbb{R}^d}\|y\|^4dQ(y),$$
which, combined with the above estimate yields
$$E\Big(\int_{C_1'}\|y\|dQ(y)\Big)^4\leq \frac 1 {n^4} \int_{\mathbb{R}^d}\|y\|^4dQ(y).$$
From this bound, (\ref{bound1}) and Schwarz's inequality we obtain
$$E(Z-Z')_+^2\leq \frac 8 {n^2} \Big(E(\|X_1-X_1'\|^2 \|X_1\|^2)+(E\|X_1-X_1'\|^4)^{1/2} \Big({\textstyle \int_{\mathbb{R}^d}\|y\|^4dQ(y) } \Big)^{1/2} \Big).$$
This and the Efron-Stein inequality for variances complete the proof. 

\hfill $\Box$

\bigskip
Theorem \ref{EfronSteindimd} provides a simple bound with explicit constants for the variance of $\mathcal{W}_2^2(P_n,Q)$
and implies tightness of $\sqrt{n}(\mathcal{W}_2^2(P_n,Q)-E(\mathcal{W}_2^2(P_n,Q)))$ with the only requirement of finite 
fourth moments and a density for $Q$. Next, we present a different application of the Efron-Stein inequality that will result in an
approximation bound from which a CLT can be concluded.

\begin{Theorem}\label{EfronSteinLin}
Assume that $P$ and $Q$ satisfy (\ref{regularQ}) and have finite 
moments of order $4+\delta$ for some $\delta>0$. Write $\varphi_0$ for the optimal transportation potential from
$P$ to $Q$. If
$$R_n=\mathcal{W}_2^2(P_n,Q)-\int_{\mathbb{R}^d} (\|x\|^2-2\varphi_0(x)) dP_n(x),$$
then
$$n \mbox{\em Var}(R_n)\to 0$$
as $n\to \infty$.
\end{Theorem}

\medskip
\noindent
\textbf{Proof.} We will argue as in the proof of Theorem \ref{EfronSteindimd}.
We write $\psi_0=\varphi_0^*$ for the optimal transportation potential from $Q$ to $P$. Without
loss of generality we can assume that $X_i=\nabla \psi_0 (U_i)$, $i=1, \ldots,n$, $X_1'=\nabla \psi_0 (U_1')$, with $U_1,\ldots,U_n,U_1'$
i.i.d. r.v.'s with law $Q$. We note that, with probability one, $\mathcal{W}_2(P_n,P)\to 0$ and we can apply Theorem~\ref{convergenceofpotentials}.
Hence, if write $\psi_n$ for the suitable centered optimal transportation potentials from $Q$ to $P_n$ that satisfy $\psi_n\to\psi_0$ $Q$-a.s.,
and $\varphi_n=\psi_n^*$, then 
\begin{equation}\label{convergence1}
\varphi_n(\nabla \psi_0(x))\to\varphi_0(\nabla \psi_0(x))
\end{equation}
for $Q$ almost every $x$.

Next, we write $P_n'$ for the empirical measure on $X_1',X_2,\ldots,X_n$ and 
$$R_n'=\mathcal{W}_2^2(P'_n,Q)-\int_{\mathbb{R}^d} (\|x\|^2-2\varphi_0(x)) dP'_n(x).$$ 
Now, the Efron-Stein inequality (\ref{EFsimple}) implies that it suffices to show that
\begin{equation}\label{targetES}
n^2 E(R_n-R_n')_+^2\to 0 \quad \mbox{ as } n \to \infty.
\end{equation}

We show first that $n(R_n-R_n')_+\to 0$ a.s.. We write $\psi_n'$ for the optimal transportation potential from $Q$ to $P_n$ and $\varphi_n'=(\psi_n')^*$.

We note that 
$$\mathcal{W}_2^2(P_n,Q)=\int_{\mathbb{R}^d} (\|x\|^2-2\varphi_n(x)) dP_n(x) +\int_{\mathbb{R}^d} (\|y\|^2-2 \psi_n(y)) dQ(y)$$
and similarly for $\mathcal{W}_2^2(P_n',Q)$ replacing $(\varphi_n,\psi_n)$ with $(\varphi_n',\psi_n')$. Also, by optimality,
$$\mathcal{W}_2^2(P_n',Q)\geq\int_{\mathbb{R}^d} (\|x\|^2-2\varphi_n(x)) dP_n'(x) +\int_{\mathbb{R}^d} (\|y\|^2-2 \psi_n(y)) dQ(y).$$
Hence,
\begin{eqnarray*}
R_n-R_n'&\leq & 2 \int_{\mathbb{R}^d} (\varphi_0(x)-\varphi_n(x)) dP_n(x)-2  \int_{\mathbb{R}^d} (\varphi_0(x)-\varphi_n(x)) dP_n'(x)\\
&=& \frac 2 n \big[(\varphi_0(X_1)-\varphi_n(X_1))-(\varphi_0(X_1')-\varphi_n(X_1'))\big]\\
&=& \frac 2 n \big[(\varphi_0(\nabla \psi_0(U_1))-\varphi_n(\nabla \psi_0(U_1)))-(\varphi_0(\nabla \psi_0(U_1'))-\varphi_n(\nabla \psi_0(U_1')))\big].
\end{eqnarray*}
Combining this bound with (\ref{convergence1}) we conclude that $n(R_n-R_n')_+\to 0$ a.s., as claimed. To complete the proof it suffices 
to show that $n^2 (R_n-R_n')_+^2$ is uniformly integrable. Since
$$n (R_n-R_n')=n\big(\mathcal{W}_2^2(P_n,Q)-\mathcal{W}_2^2(P_n',Q) \big)- \big((\|X_1\|^2-2\varphi_0(X_1))-(\|X_1'\|^2-2\varphi_0(X_1')) \big),$$
and $(\|X_1\|^2-2\varphi_0(X_1))$ and $(\|X_1'\|^2-2\varphi_0(X_1'))$ have finite second moment (recall Theorem \ref{secondmomentth}), this will
follow if we prove that $n^2\big(\mathcal{W}_2^2(P_n,Q)-\mathcal{W}_2^2(P_n',Q) \big)_+^2$ is uniformly integrable. For this last goal we write 
$Z=\mathcal{W}_2^2(P_n,Q)$, $Z'=\mathcal{W}_2^2(P_n',Q)$ and recall
from the proof of Theorem \ref{EfronSteindimd} that 
$$n (Z_n-Z_n')_+\leq \|X_1-X_1'\| \Big(\big(\|X_1\|+\|X_1'\|\big)+2n\int_{C_1'}\|y\|dQ(y)\Big),$$
keeping the notation there for $C_1'$. Since $X_1,X_1'$ have finite fourth moment, we only need to prove that
$\Big(n\|X_1-X_1'\|\int_{C_1'}\|y\|dQ(y)\Big)^2$ is uniformly integrable.
To check this we argue as above to see that
$$\Big(\int_{C_1'}\|y\|dQ(y)\Big)^{4+\delta}\leq
\frac 1  {n^{3+\delta}}\Big(\int_{C_1'}\|y\|^{4+\delta}dQ(y)\Big)$$
and, as a consequence, 
$$E\Big(n\int_{C_1'}\|y\|dQ(y)\Big)^{4+\delta} \leq n E\Big(\int_{C_1'}\|y\|^{4+\delta}dQ(y)\Big)=
\int_{\mathbb{R}^d}\|y\|^{4+\delta}dQ(y)<\infty.$$
Finally, we use Schwarz's inequality to see that
$$E\Big(n\|X_1-X_1'\|\int_{C_1'}\|y\|dQ(y)\Big)^{2+ \frac \delta 2}\leq \Big(E \|X_1-X_1'\|^{4+\delta}\Big)^{\frac 1 2} \Big(
\int_{\mathbb{R}^d}\|y\|^{4+\delta}dQ(y)\Big)^{\frac 1 2}.$$
This entails that $\Big(n\|X_1-X_1'\|\int_{C_1'}\|y\|dQ(y)\Big)^2$ is uniformly integrable and completes the proof.
\hfill $\Box$

We consider next a version of the variance bounds in Theorems \ref{EfronSteindimd} and \ref{EfronSteinLin} suited to the two-sample 
empirical transportation cost. Thus, we assume that $X_1,\ldots,X_n$ are i.i.d. r.v.'s with law $P$,
$Y_1,\ldots,Y_m$ are i.i.d. r.v.'s with law $Q$, independent of the $X_i$'s, $P_n$ denotes the empirical
measure on the $X_i$'s and $Q_m$ the empirical measure on the $Y_j$'s.

\begin{Theorem}\label{ESTS}
If $P$ and $Q$ have densities and finite fourth moments then
$$\mbox{\em Var}(\mathcal{W}_2^2(P_n,Q_m))\leq \frac{C(P,Q)}{n}+\frac{C(Q,P)}{m},$$
where $C(P,Q)$ is defined as in Theorem \ref{EfronSteindimd}.

If $P$ and $Q$ satisfy (\ref{regularQ}) and have finite moments of order  $4+\delta$ for some $\delta>0$, $n\to\infty$, $m\to \infty$,
$\frac n {n+m} \to \lambda \in (0,1)$ and set

$$R_{n,m}=\mathcal{W}_2^2(P_n,Q_m)-\int_{\mathbb{R}^d} (\|x\|^2-2\varphi_0(x))dP_n(x)-\int_{\mathbb{R}^d} (\|y\|^2-2\psi_0(y))dQ_m(y),$$
then
$$\frac{nm}{n+m} \mbox{\em Var}(R_{n,m})\to 0.$$
\end{Theorem}

\medskip
\noindent \textbf{Proof.} We note first that, as a function of $X_1,\ldots,X_n,Y_1,\ldots,Y_m$, 
$\mathcal{W}_2^2(P_n,Q_m)$ is symmetric in its first $n$ variables, as well as in its last $m$. Hence, using the Efron-Stein inequality
we see that
$$\mbox{Var}(\mathcal{W}_2^2(P_n,Q_m))\leq n E(Z-Z')_+^2+mE(Z-Z'')_+^2,$$
where $Z=\mathcal{W}_2^2(P_n,Q_m)$, $Z'=\mathcal{W}_2^2(P_n',Q_m)$, $Z''=\mathcal{W}_2^2(P_n,Q_m')$,
$P_n'$ is the empirical measure on $X_1',X_2,\ldots,X_n$, $Q_m'$ is the empirical measure on $Y_1',Y_2,\ldots,Y_m$
and $X_1',Y_1'$ are independent r.v.'s, independent of the $X_i$'s and $Y_j$'s, with $X_1'$ having law $P$ and $Y_1'$ with
law $Q$. To bound $E(Z-Z')_+^2$ we write $\pi$ (resp. $\pi'$) for the optimal transportation plan from $P_n$
to $Q_m$ (resp. from $P_n'$ to $Q_m$). We write also $\pi_{i,j}$ for the probability that $\pi$ assigns to the 
pair $(X_i,Y_j)$, and similarly for $\pi_{i,j}'$, $c_{i,j}=\|X_i-Y_j\|^2$ and $c_{i,j}'$ for the
costs associated to the data $X_1',X_2,\ldots,X_n,Y_1,\ldots,Y_m$. Then $Z'=\sum_{i=1}^n\sum_{j=1}^m c'_{i,j}\pi'_{i,j}$
and $Z\leq \sum_{i=1}^n\sum_{j=1}^m c_{i,j}\pi'_{i,j}$. Hence, noting that $c_{i,j}=c'_{i,j}$ for $i\geq 2$ we see that
$$Z-Z'\leq \sum_{j=1}^m \pi'_{1,j}(c_{1,j}-c'_{1,j})\leq \|X_1-X_1'\| \sum_{j=1}^m \pi'_{1,j}(\|X_1\|+\|X_1\|'+2\|Y_j\|).$$
Since $\sum_{j=1}^m \pi'_{1,j}=\frac 1 n$ we obtain that
$$Z-Z'\leq \|X_1-X_1'\| \big( \frac 1 n (\|X_1\|+\|X_1\|')+2\sum_{j=1}^m \pi'_{1,j}\|Y_j\| \big).$$
From this point we can argue as in the proof of Theorem \ref{EfronSteindimd} to conclude that $E(Z-Z')_+^2\leq \frac{C(P,Q)}{n^2}$.
We note that, again in this setup, we have by exchangeability
$$E\Big(\sum_{j=1}^m \pi'_{1,j}\|Y_j\|^4\Big)=\frac 1 n E\Big(\sum_{i=1}^n\sum_{j=1}^m \pi'_{i,j}\|Y_j\|^4\Big)=
\frac 1 n E\Big(\frac 1 m \sum_{j=1}^m \|Y_j\|^4\Big)=\frac 1 n E \|Y_1\|^4.$$
Similarly, we see that $E(Z-Z'')_+^2\leq \frac{C(Q,P)}{m^2}$ and this proves the first claim.

For the second claim we argue as in the proof of Theorem \ref{EfronSteinLin}. We keep the notation $P_n'$, $Q_m'$ as above and
set 
$$R'_{n,m}=\mathcal{W}_2^2(P_n',Q_m)-\int_{\mathbb{R}^d} (\|x\|^2-2\varphi_0(x))dP'_n(x)-\int_{\mathbb{R}^d} (\|y\|^2-2\psi_0(y))dQ_m(y),$$
$$R''_{n,m}=\mathcal{W}_2^2(P_n,Q_m')-\int_{\mathbb{R}^d} (\|x\|^2-2\varphi_0(x))dP_n(x)-\int_{\mathbb{R}^d} (\|y\|^2-2\psi_0(y))dQ'_m(y).$$
Again, the Efron-Stein inequality shows that it suffices to prove that
$n^2E(R_{n,m}-R'_{n,m})_+^2\to 0$ and $m^2E(R_{n,m}-R''_{n,m})_+^2\to 0$. We prove the first of these two claims, the other 
following by symmetry. We write $\varphi_n$ for the optimal transportation potential from $P_n$ to $Q_m$ and $\psi_n=\varphi_n^*$. We note that
Theorem \ref{convergenceofpotentials} ensures that we can center the $\phi_n$'s to ensure that $\varphi_n\to\varphi_0$ $P$-a.s.. Also, as above,
$$\mathcal{W}_2^2(P_n,Q_m)=\int_{\mathbb{R}^d} (\|x\|^2-2\varphi_n(x)) dP_n(x) +\int_{\mathbb{R}^d} (\|y\|^2-2 \psi_n(y)) dQ_m(y),$$
while
$$\mathcal{W}_2^2(P_n',Q_m)\geq\int_{\mathbb{R}^d} (\|x\|^2-2\varphi_n(x)) dP_n'(x) +\int_{\mathbb{R}^d} (\|y\|^2-2 \psi_n(y)) dQ_m(y).$$
From this we see that
\begin{eqnarray*}
R_{n,m}-R_{n,m}'&\leq & 2 \int_{\mathbb{R}^d} (\varphi_0(x)-\varphi_n(x)) dP_n(x)-2  \int_{\mathbb{R}^d} (\varphi_0(x)-\varphi_n(x)) dP_n'(x)\\
&=& \frac 2 n \big[(\varphi_0(X_1)-\varphi_n(X_1))-(\varphi_0(X_1')-\varphi_n(X_1'))\big]
\end{eqnarray*}
and this shows that $n(R_{n,m}-R_{n,m}')_+\to 0$ a.s.. Arguing as in the proof of Theorem \ref{EfronSteinLin} we can check that 
$n^2(R_{n,m}-R_{n,m}')_+$ is uniformly integrable. Hence, we conclude that $n^2E(R_{n,m}-R'_{n,m})_+^2\to 0$ and complete the proof.
\hfill $\Box$

\section{CLTs for empirical transportation cost} \label{CLTsequel}
As a direct consequence of the approximation bounds in Theorems \ref{EfronSteinLin} and \ref{ESTS} we arrive to the main results in this paper, namely, central limit theorems for the empirical transportation cost and the optimal matching cost.

\begin{Theorem}[Central Limit Theorem for empirical quadratic transportation cost]\label{CLTdimd}
Assume $P$ and $Q$ are probabilities on $\mathbb{R}^d$ that satisfy (\ref{regularQ}) and have finite moments of order $4+\delta$
for some $\delta>0$. If $X_1,\ldots,X_n$ are i.i.d. r.v.'s with law $P$ and
$P_n$ denotes the empirical measure on $X_1,\ldots,X_n$ then
$$n\mbox{\em Var}(\mathcal{W}_2^2(P_n,Q))\to \sigma^2(P,Q):=\int_{\mathbb{R}^d} (\|x\|^2-2\varphi_0(x))^2 dP(x)-
\Big(\int_{\mathbb{R}^d} (\|x\|^2- 2 \varphi_0(x)) dP(x)\Big)^2$$
and
$$\sqrt{n}\big(\mathcal{W}_2^2(P_n,Q)-E\mathcal{W}_2^2(P_n,Q))\underset w \rightarrow N(0,\sigma^2(P,Q))$$
as $n\to\infty$, where $\varphi_0$ denotes an optimal transportation potential from $P$ to $Q$.

Furthermore, if $Y_1,\ldots,Y_m$ are i.i.d. r.v.'s with law $Q$, independent of the $X_i$'s, 
$Q_m$ denotes the empirical measure on $Y_1,\ldots,Y_m$ and $n\to\infty$, $m\to\infty$ with
$\frac n {n+m}\to \lambda \in (0,1)$, then
$${\textstyle\frac {nm}{n+m}}\mbox{\em Var}(\mathcal{W}_2^2(P_n,Q_m))\to (1-\lambda)\sigma^2(P,Q)+\lambda \sigma^2(Q,P)$$
and
$${\textstyle\sqrt{\frac {nm}{n+m}}}\big(\mathcal{W}_2^2(P_n,Q_m)-E\mathcal{W}_2^2(P_n,Q_m))\underset w \rightarrow 
N(0,(1-\lambda)\sigma^2(P,Q)+\lambda \sigma^2(Q,P)).$$
\end{Theorem}

\medskip
We believe that the assumptions of moments with order $4+\delta$ is a technical condition that could be weakened to moments of order 4 only. Yet, for the proof, this condition is mandatory.

To end this Section we provide an additional CLT for $\mathcal{W}_2^2(P_n,Q)$ which does not require smoothness on $P$, but only
on $Q$. Now a finite fourth moment for $Q$ will suffice, but $P$ will be assumed to have finite support. The proof will use
the following special form for the quadratic transportation cost to a finitely supported probability.

\begin{Proposition}\label{specialform}
Assume $P$ has finite support, $\{x_1,\ldots,x_k\}\subset \mathbb{R}^d$, with $P\{x_i\}=p_i $, $i=1,\ldots,k$
and $Q$ is a Borel probability on $\mathbb{R}^d$ with finite second moment then
$$\mathcal{W}_2^2(P,Q)=\int_{\mathbb{R}^d} \|x\|^2dP(y) + \int_{\mathbb{R}^d} \|y\|^2dQ(y)-2 \min_{z\in \mathbb{R}^k}V(z),$$
where $V$ is the convex function
\begin{equation}\label{cV} V(z_1,\ldots,z_k)=\sum_{i=1}^k p_iz_i+E\max_{1\leq j\leq k} \big(x_j\cdot Y-z_j\big), \end{equation}
and $Y$ is a random vector with distribution $Q$.

If $Q\ll \ell^d$, the $d$-dimensional Lebesgue measure, then $V$ is differentiable and 
$$\nabla V(z)=(p_1,\ldots,p_k)-(Q(A_1(z)),\ldots,Q(A_k(z))),$$
where
$$A_j(z)=\big\{y\in\mathbb{R}^d:\, (x_j\cdot y-z_j)> \max_{i\ne j}(x_i\cdot y-z_i) \big\}, \quad j=1,\ldots,n.$$
Finally, if $Q$ satisfies (\ref{regularQ}) then 
$z$ minimizes $V$ if and only $\nabla V(z)=0$ and there is a unique $z$ such that $\nabla V(z)=0$, $z_i+\frac {\|x_i\|^2}{2}\geq 0$, $i=1,\ldots,k$ and
$\sum_{i=1}^k p_i (z_i+\frac {\|x_i\|^2}{2})=\max_{1\leq i\leq k}\|x_i\|^2+\int_{\mathbb{R}^d}\|y\|^2 dQ(y)$.

\end{Proposition}

\bigskip
\noindent
\textbf{Proof.} From duality theory for optimal transportation we know that
$$\mathcal{W}_2^2(P,Q)=\sum_{i=1}^k p_i \|x_i\|^2 + \int_{\mathbb{R}^d} \|y\|^2dQ(y)-2\min_{(z,\psi)\in \Phi}
\Big[\sum_{i=1}^k p_iz_i +\int \psi(y) dQ(y) \Big],$$
where $\Phi$ is the class of pairs $(z,\psi)$ such that $z\in\mathbb{R}^k$, $\psi\in L_1(Q)$ and
$$x_j\cdot y\leq z_j+\psi(y),\quad  1\leq j\leq k, y\in\mathbb{R}^k.$$
Since
$$\psi(y)\geq \tilde{\psi}(y):=\max_{1\leq j\leq k} (x_j\cdot y-z_j)$$
and $(z,\tilde{\psi})\in \Phi$ we see that
$$\min_{(z,\psi)\in \Phi}\Big[\sum_{i=1}^k p_iz_i +\int \psi(y) dQ(y) \Big]=\min_{z\in\mathbb{R}^d} V(z)$$
with $V$ as in the statement \eqref{cV}, which is obviously convex. Let us fix now $z\in\mathbb{R}^k$, set $\psi(y)=
\max_{1\leq j\leq k} (x_j\cdot y-z_j)$ and consider
$\tilde{z}_j=\sup_{y\in\mathbb{R}^d} (x_j\cdot y-\psi(y))$. 
Since $z_j\geq x_j\cdot y-u(y)$ for all $y$ we have $\tilde{z}_j\leq z_j$, $j=1,\ldots,n$. Let us 
now set  $\tilde{\psi}(y)=\max_{1\leq j\leq n}(x_j\cdot y-\tilde{z}_j)$. Then we have 
$\tilde{\psi}(y)=\max_{1\leq j\leq n}(x_j\cdot y-\tilde{z}_j)\geq \max_{1\leq j\leq n}(x_j\cdot y-{z}_j)=\psi(y)$. On the other hand,
$\tilde{z}_j+\psi(y)\geq x_j\cdot y$ for all $j$ and $y$ implies  $\psi(y)\geq \max_{1\leq j\leq n} (x_j\cdot y - \tilde{z}_j)= \tilde{\psi}(y)$.
Hence, $\tilde{\psi}=\psi$ and $V(\tilde{z}_1,\ldots,\tilde{z}_k)\leq V({z}_1,\ldots,{z}_k)$. If $p_i>0$ then the last inequality is
strict unless $\tilde{z}_i=z_i$. 

From this point we assume that $Q$ has a density. Then a minimizing pair 
$(z,\psi)$ in $\Phi$ must satisfy $z_j=\sup_{y\in\mathbb{R}^d} (x_j\cdot y-\psi(y))$ and
$\nabla \psi$ is the optimal transportation map from $Q$ to $P$. Since, on the other hand,
$\psi(y)=\max_{1\leq j\leq n}(x_j\cdot y-z_j)$ we see that $\nabla \psi(y)=x_j$ if $y\in A_j(z)$ and 
the condition $Q(A_j(z))=p_j$, $j=1,\ldots,k$ is necessary and sufficient for $z$ to be a minimizer of $V$.

If $Q$ satisfies (\ref{regularQ}) then the polyhedral sets that
are mapped by $\nabla \psi$ onto the $x_i$'s are uniquely determined up to differences in the boundaries,
which entails that any two minimizers $u,\tilde{\psi}$ satisfy $\tilde{\psi}=\psi+L$ for some constant $L$. Consequently,
two minimizers, $z, \tilde{z}$ of $V$ must satisfy $\tilde{z}_i=z_i-L$, $i=1,\ldots,k$.

For the claims about the differentiability of $V$ it suffices to focus on 
$$\tilde{V}(z)=E\max_{1\leq j\leq k} \big(x_j\cdot Y-z_j\big)$$ and note that
\begin{eqnarray*}
\lefteqn{\tilde{V}(z+h)-\tilde{V}(z)-\sum_{j=1}^k h_j Q(A_j(z))}\hspace*{3cm}\\
&=&\sum_{j=1}^k E\Big[\Big( \max_{1\leq i\leq k} \big(x_i\cdot Y-(z_{i}+h_i)\big)- 
\big(x_j\cdot Y-(z_{j}+h_j)\big) \Big)I_{A_j(z)}(Y)\Big].
\end{eqnarray*}
It is easy to check that $0\leq \Big( \max_{1\leq i\leq k} \big(x_i\cdot Y-(z_{i}+h_i)\big)- 
\big(x_j\cdot Y-(z_{j}+h_j)\big) \Big)I_{A_j(z)}(Y)\leq 2 \max_{1\leq j\leq k}|h_j|$, while, as $h\to 0$,
$\big(\max_{1\leq i\leq k} \big(x_i\cdot Y-(z_{i}+h_i)\big)- 
\big(x_j\cdot Y-(z_{j}+h_j)\big) I_{A_j(z)}(Y)$ eventually vanishes (except, possibly, if 
$Y$ belongs to the boundary of $A_j(z)$). Then, from dominated convergence we conclude that
$$\frac{\tilde{V}(z+h)-\tilde{V}(z)-\sum_{j=1}^k h_j Q(A_j(z))}{\|h\|}\to 0$$
as $\|h\|\to 0$, proving that $\tilde{V}$, and therefore, $V$ are differentiable.
Obviously, the condition $\nabla V(z)=0$ is exactly the necessary and sufficient condition 
for $z$ to be a minimizer of $V$ shown above.

Finally, let us fix $z\in\mathbb{R}^k$ and write $\psi(y)=\max_{1\leq j\leq k} (x_j\cdot y-z_j)$.
Since
$$\psi(y)+\frac{\|y\|^2}{2}\geq \frac{\|x_j+y\|^2}{2}-z_j-\frac{\|x_j\|^2}{2}\geq -z_j-\frac{\|x_j\|^2}{2}$$
we see that $a:=\inf_{y\mathbb{R}^d} \big(\psi(y)+\frac{\|y\|^2}{2}\big)$ $a$ is finite.
As noted above, $V$ remains unchanged if we replace $(z_1,\ldots,z_k)$ by $(z_1+a,\ldots,z_k+a)$ and $\psi(y)$ becomes $\psi(y)-a$.
As a consequence, in the minimisation of $V$ it suffices to consider points $(z_1,\ldots,z_k)$ such that
\begin{equation}\label{inf0}
\inf_{y\in\mathbb{R}^d} \Big(\psi(y)+\frac{\|y\|^2}2\Big)=0.
\end{equation}
Let us assume that (\ref{inf0}) holds and
consider $\tilde{z}_j=\sup_{y\in\mathbb{R}^d} (x_j\cdot y-\psi(y))$. As above, we have
$\tilde{z}_j\leq z_j$, $j=1,\ldots,n$, $\tilde{\psi}(y)=\max_{1\leq j\leq k}(x_j\cdot y-\tilde{z}_j)=\psi(y)$ 
and $V(\tilde{z}_1,\ldots,\tilde{z}_k)\leq V({z}_1,\ldots,{z}_k)$.
We observe now that
\begin{eqnarray*}
\tilde{z}_j+\frac{\|x_j\|^2}{2}&=& \sup_{y\in\mathbb{R}^d} \Big(x_j\cdot y +\frac{\|x_j\|^2}2-\psi(y)\Big)\\
&\geq &
\sup_{y\in\mathbb{R}^d} \Big(-\frac{\|y\|^2}2-\psi(y)\Big)=-\inf_{y\in\mathbb{R}^r} \Big(\psi(y)+\frac{\|y\|^2}2\Big)=0.
\end{eqnarray*}
On the other hand,
$$V(\tilde{z}_1,\ldots,\tilde{z}_k)+\frac 1 2\sum_{i=1}^k p_i\|x_i\|^2+\frac 1 2 \int_{\mathbb{R}^d}\|y\|^2 dQ(y)=
\sum_{i=1}^k p_i\Big( \tilde{z}_i+\frac {\|x_i\|^2}2 \Big)+\int_{\mathbb{R}^d} \psi(y)+\frac{\|y\|^2}2 dQ(y),
$$
which, by (\ref{inf0}), implies that 
$$\sum_{i=1}^k p_i\Big( \tilde{z}_i+\frac {\|x_i\|^2}2 \Big)\leq 
V(\tilde{z}_1,\ldots,\tilde{z}_k)+\frac 1 2\sum_{i=1}^k p_i\|x_i\|^2+\frac 1 2 \int_{\mathbb{R}^d}\|y\|^2 dQ(y).$$
Nonnegativity of $\mathcal{W}_2^2(P,Q)$ shows that $\min_{z\in\mathbb{R}^k} V(z)\leq \frac 1 2 \Big( \sum_{i=1}^k p_i\|x_i\|^2+
\int_{\mathbb{R}^d} \|y\|^2 dQ(y) \Big)$. Hence, there exists a minimizer of $V$ that satisfies $z_i+\frac{\|x_i\|^2}{2}\geq 0$ and 
$\sum_{i=1}^k p_i( {z}_i+\frac {\|x_i\|^2}2 )\leq M:=\max_{1\leq i\leq k}\|x_i\|^2$ $+\int_{\mathbb{R}^d}\|y\|^2 dQ(y)$. 
Adding a constant, if necessary, we see that there is a unique minimizer of $V$ that satisfies $z_i+\frac{\|x_i\|^2}{2}\geq 0$ and 
$\sum_{i=1}^k p_i( {z}_i+\frac {\|x_i\|^2}2 )= M$.

\hfill \quad $\Box$

\bigskip
We note that the minimizing $z=(z_1,\ldots,z_k)$ in Proposition~\ref{specialform} 
satisfy $z_i=\varphi_0(x_i)$, $i=1,\ldots,k$ with $\varphi_0=\psi_0^*$ and $\psi_0$ the optimal transportation
potential from $Q$ to $P$ (which is unique up to the addition of a constant by Theorem \ref{uniquepotential} 
under (\ref{regularQ}). Hence, we see that the optimal transportation potential from $P$ to $Q$ is
also unique (up to the addition of a constant) in this setup.

We can prove now the announced CLT for $\mathcal{W}_2^2(P_n,Q)$ when $P$ is finitely supported.

\begin{Theorem}\label{CLTfinitesupport}
If $P$ has a finite support, and moreover if $Q$ satisfies (\ref{regularQ}) and has a finite fourth moment. If
$X_1,\ldots,X_n$ are i.i.d. r.v.'s with law $P$ and $P_n$ denotes the empirical measure on $X_1,\ldots,X_n$, then
$$\sqrt{n}\big(\mathcal{W}_2^2(P_n,Q)-\mathcal{W}_2^2(P,Q)\big)\underset w \rightarrow N(0,\sigma^2(P,Q))$$
as $n\to\infty$, where $\sigma^2(P,Q)$ is as in Theorem \ref{CLTfinitesupport}.
\end{Theorem}

\medskip
\noindent \textbf{Proof.} 
We assume that $P$ is as in Proposition~\ref{specialform}. We can write
$\mathcal{W}_2^2(P,Q)=\max_{z\in C_M} M(z)$ with
$$M(z)=\sum_{j=1}^k p_j \|x_j\|^2+\int_{\mathbb{R}^d} \|y\|^2 dQ(y)-2 \sum_{j=1}^kp_j z_j-2\tilde{V}(z),$$
$\tilde{V}(z)=E\max_{1\leq j\leq k} \big(x_j\cdot Y-z_j\big)$ and $C_M=\{z\in\mathbb{R}^d:\, z_i+\frac{\|x_i\|^2}{2}\geq 0, i=1,\ldots,k;\, 
\sum_{i=1}^k p_i( {z}_i+\frac {\|x_i\|^2}2 )= M\}$. Similarly, $\mathcal{W}_2^2(P_n,Q)=
\max_{z\in C_M} M_n(z)$, where $M_n$ is obtained replacing the $p_j$'s by the empirical frequencies,
$p_{n,j}$'s. We write $z_n$ and $z_0$ for the unique maximizers of $M_n$ and $M$, respectively, given by in Theorem \ref{specialform}. 
By the Central Limit Theorem in $\mathbb{R}^k$ we have $U_n:=[\sqrt{n}(p_{n,j}- p_j)]_{1\leq j\leq k}\underset w \to U$ with $U$ a centered
Gaussian random vector with covariance matrix $\Sigma=[\sigma_{i,j}]_{1\leq i,j\leq k}$, $\sigma_{i,i}=p_i (1-p_i)$, $\sigma_{i,j}=-p_i p_j$, $i\neq j$.
Without loss of generality we can assume that $U_n\to U$ a.s.. Note that, in particular, 
\begin{equation}\label{eq2}
M_n(z)-M(z)=\frac 1 {\sqrt{n}}\sum_{i=1}^k U_{n,i}(\|x_i\|^2-2z_{i}).
\end{equation}
On the other hand, the choice of $z_n$ guarantees that it is a bounded sequence. Assume that, through a subsequence, $z_n\to\hat{z}$.
Then $M_n(z_n)\to M(\hat{z})$ (here we are using the continuity of $\tilde{V}$. For any fixed $z$ we have
$M(z)=\lim_{n\to\infty} M_n(z)\leq\lim_{n\to\infty} M_n(z_n)=M(\hat{z})$. Hence, $\hat{z}$ is a maximizer of $M$. But obviously
$\hat{z}_i+\frac{\|x_i\|^2} 2\geq 0$ and $\sum_{i=1}^kp_i(\hat{z}_i+\frac{\|x_i\|^2} 2)=M$. Hence,
by uniqueness, we must have $\hat{z}=z_0$, that is, $z_n\to z_0$ a.s.. From this fact we see that
\begin{eqnarray}\nonumber
\sqrt{n}(\mathcal{W}_2^2(P,Q)-\mathcal{W}_2^2(P,Q))&=&\sqrt{n}(M_n(z_n)-M(z_0))\\
&=&\sqrt{n}(M_n(z_n)-M(z_n))+\sqrt{n}(M(z_n)-M(z_0)).\label{eq1}
\end{eqnarray}
Now, by optimality we see that $\sqrt{n}(M_n(z_0)-M(z_0))-\sqrt{n}(M_n(z_n)-M(z_n))\leq \sqrt{n}(M(z_n)-M(z_0))\leq 0$. Also, from 
(\ref{eq2}) we see that $\sqrt{n}(M_n(z_n)-M(z_n))\to \sum_{i=1}^k U_{i}(\|x_i\|^2-2z_{0,i})$, 
$\sqrt{n}(M_n(z_0)-M(z_0))\to \sum_{i=1}^k U_{i}(\|x_i\|^2-2z_{0,i})$ a.s.. As a consequence, $\sqrt{n}(M(z_n)-M(z_0))\to 0$ a.s. which,
together with (\ref{eq1}), shows that
$$\sqrt{n}(\mathcal{W}_2^2(P_n,Q)-\mathcal{W}_2^2(P,Q))\underset w \to \sum_{i=1}^k U_{i}(\|x_i\|^2-2z_{0,i}).$$
A simple computation shows that the right hand side in this last display is a centered Gaussian random variable
with variance $\sigma^2(P,Q)$ as in Theorem \ref{EfronSteindimd}.

\hfill $\Box$

\begin{Remark}{\em
We note that, provided $Q$ has a finite moment of order $4+\delta$ for some $\delta>0$, the linearization bound 
in Theorem \ref{EfronSteinLin} can be adapted to cover this setup and conclude that
$$n\mbox{Var}(\mathcal{W}_2^2(P_n,Q))\to \sigma^2(P,Q)$$
and
$$\sqrt{n}\big(\mathcal{W}_2^2(P_n,Q)-E(\mathcal{W}_2^2(P_n,Q))\big)\underset w \rightarrow N(0,\sigma^2(P,Q)).$$
On the other hand, the centering constants $E(\mathcal{W}_2^2(P_n,Q))$ in Theorem \ref{CLTdimd} cannot be replaced 
in general by $\mathcal{W}_2^2(P,Q)$. As an example, consider the case when $P=Q$ is the uniform distribution on 
the $d$-dimensional unit cube. In this case 
Theorem \ref{CLTdimd} yields that
$$\sqrt{n}\big( \mathcal{W}_2^2(P_n,Q)-E \mathcal{W}_2^2(P_n,Q) \big)\to 0$$
in probability.
% On the other hand $E \mathcal{W}_2^2(P_n,Q)$ is of order $n^{-1/d}$ if $d\geq 3$ (see \cite{Talagrand}) and
%we cannot have $\sqrt{n} \mathcal{W}_2^2(P_n,Q)\to 0$ (otherwise we would conclude that $E \mathcal{W}_2^2(P_n,Q)=o(n^{-1/2})$).
 On the other hand $E\mathcal{W}_2^2(P_n,Q)$ is of order $n^{-2/d}$ if $d\geq 5$ (see Theorem 1 and subsequent comments in \cite{Fournierguillin}) and we cannot have $\sqrt\mathcal{W}_2^2(P_n,Q)\to 0$ (otherwise we would conclude that $E\mathcal{W}_2^2(P_n,Q)=o(n^{-1/2})$.}
 \end{Remark}
 
 To conclude, we would like to add two final comments. First, we note that in the case $P=Q$ Theorem 1 in \cite{Fournierguillin} yields that (provided $d\geq 5$ and assuming that $P$ has finite moment of order $q>\frac{2d}{d-2}$)

$$n^{d/2}(\mathcal{W}_2^2(P_n,Q)-E\mathcal{W}_2^2(P_n,Q))$$

is stochastically bounded. In this setup, assuming  $P$ has finite moment of order 4 (and a density) we see that 

$$\sqrt{n}(\mathcal{W}_2^2(P_n,Q)-E\mathcal{W}_2^2(P_n,Q))$$

is stochastically bounded. Under slightly stronger assumptions, Theorem 4.1, shows that 

$$\sqrt{n}(\mathcal{W}_2^2(P_n,Q)-E \mathcal{W}_2^2(P_n,Q))\to 0
$$
in probability. It would be of great interest to investigate whether a nontrivial CLT holds in this setup at a different rate.\\ \indent In the one-dimensional case the problem was considered in \cite{dBGU05}, proving weak convergence to some non degenerate and non Gaussian limit law. This case provides some indication that the case $P=Q$ is, essentially, of a different nature and that a nontrivial CLT in that case cannot be obtained with the techniques used in this paper.

\hfill $\Box$


\begin{thebibliography}{99}


\bibitem{AjtaiKomlosTusnady} Ajtai , M., Koml\'os , J. and Tusn\'ady , G. (1984). On optimal matchings. \textit{Combinatorica}, \textbf{4}
259--264.

\bibitem{ambrosio2016pde} Ambrosio, L., Stra, F.  and Trevisan, D. (2016).  {A PDE approach to a 2-dimensional matching problem}. \textit{arXiv preprint arXiv:1611.04960}.

\bibitem{delBarrioMatran2013a} del Barrio, E. and Matr\'an, C. (2013). Rates of convergence for partial mass problems. 
\textit{Probab. Theory Relat. Fields}, \textbf{155}, 521--542.

\bibitem{dBGM99} del Barrio, E., Matr\'an, C. and Gin\'e, E. (1999). Central limit theorems for the 
Wasserstein distance between the empirical and the true distribution. \textit{Ann. Probab.}, \textbf{27}, 1009--1071.

\bibitem{dBGU05} del Barrio, E., Gin\'e, E. and Utzet, F. (2005). Asymptotics for $L_2$ functionals of the empirical 
quantile process, with applications to tests of fit based on weighted Wasserstein distances. \textit{Bernoulli},
\textbf{11}, 131--189.


\bibitem{BovkovLedoux2014} Bobkov, S. and Ledoux, M. (2016). One-dimensional empirical measures, 
order statistics and Kantorovich transport distances. \textit{To appear in Memoirs of the AMS.}

\bibitem{BoucheronLugosiMassart2013} Boucheron, S., Lugosi, G. and Massart, P. (2013). \textit{Concentration Inequalities. A Nonasymptotic
Theory of Independence.} Oxford.


\bibitem{CuestaAlbertosMatranTueroDiaz} Cuesta-Albertos, J.A., Matr\'an, C. and Tuero-D\'{\i}az, A. (1997).
Optimal transportation plans and convergence in distribution. \textit{J. Multivariate Analysis}, \textbf{60}, 72--83.


\bibitem{DobricYukich} Dobri\'c, V. and Yukich , J. E. (1995). Asymptotics for transportation cost in high dimensions. 
\textit{J. Theoret. Probab.}, \textbf{8}, 97--118.

\bibitem{GariepyEvans} Evans, L. C. and Gariepy, R. F. (1992). \textit{Measure Theory and Fine Properties of Functions.} CRC Press.

\bibitem{Fournierguillin} Fournier, N. and Guillin, A. (2015). {On the rate of convergence in Wasserstein distance of the empirical measure}, \textit{Probability Theory and Related Fields},
  \textbf{162}, {3-4}, {707--738}.

\bibitem{GangboMcCann} Gangbo, W. and McCann, R. J. (1996). The Geometry of Optimal Transportation. \textit{Acta Math.}, \textbf{177},
113--161.

\bibitem{Heinich} Heinich, H., and Lootgieter, J. C. (1996). Convergence des fonctions monotones. Comptes rendus de l'Acad\'emie des sciences. S\'erie 1, Math\'ematique, 322(9), 869-874.

\bibitem{RachevRuschendorf}  Rachev , S. T. and R\"uschendorf , L. (1998). \textit{Mass Transportation Problems. (2 Vols).}
Springer, New York.


\bibitem{RipplMunkSturm} Rippl, T., Munk, A. and Sturm, A. (2016). Limit laws of the empirical Wasserstein distance: Gaussian distributions.
\textit{J. of Multivariate Analysis}, \textbf{151}, 90--109.

\bibitem{Rockafellar} Rockafellar, R. T. (1970). \textit{Convex Analysis}. Princeton University Press.

\bibitem{RockafellarWets} Rockafellar, R. T. and Wets, R. J. (1998). \textit{Variational Analysis}. Springer.


\bibitem{SommerfeldMunk} Sommerfeld, M. and Munk, A. (2016). Inference for Empirical Wasserstein Distances on Finite Spaces.
\textit{Preprint.} https://arxiv.org/abs/1610.03287v1


\bibitem{Talagrand1992} Talagrand, M. (1992). Matching random samples in many dimensions. \textit{Ann. Appl. Probab.},
\textbf{2}, 846--856.

\bibitem{Talagrand} Talagrand , M. (1994). The transportation cost from the uniform measure to the empirical
measure in dimension $\geq 3$. \textit{Ann. Probab.}, \textbf{22}, 919--959.

\bibitem{TalagrandYukich} Talagrand , M. and Yukich , J. E. (1993). The integrability of the square exponential 
transportation cost. \textit{Ann. Appl. Probab.}, \textbf{3}, 1100--1111.

\bibitem{Villani} Villani, C. (2003). \textit{Topics in Optimal Transportation.} American Mathematical Society.

\bibitem{VillaniBig} Villani, C. (2009). \textit{Optimal Transport: Old and New. Grundlehren der Mathematischen Wissenschaften 
[Fundamental Principles of Mathematical Sciences]} \textbf{338}. Springer, Berlin.


\end{thebibliography}
\end{document}